\documentclass[11pt]{article}

\usepackage[alphabetic]{amsrefs}
\usepackage{amsmath,amssymb,amsthm,enumerate}
\usepackage[all,cmtip]{xy}
\def\today{\number\day .\number\month .\number\year}

\def \a{{\mathfrak a}}
\def\adm{{\rm adm}}
\def \al{\alpha}
\def \aut{{\rm aut}}
\def \bs{\backslash}
\def \C{{\mathbb C}}
\def \CD{{\cal D}}

\def \cusp{\mathrm{cusp}}

\def \diag{{\rm diag}}

\def\e{\emph}

\def \eps{\varepsilon}
\def \g{{\mathfrak g}}
\def \ga{\gamma}
\def \Ga{\Gamma}
\def \GL{{\rm GL}}
\def\H{{\mathbb H}}

\def \Id{{\rm Id}}
\def \Im{{\rm Im}}

\def \k{{\mathfrak k}}
\def \la{\lambda}

\def \n{{\mathfrak n}}
\def \N{{\mathbb N}}
\def \ol{\overline}
\def \p{{\mathfrak p}}
\def \PGL{{\rm PGL}}
\def \ph{\varphi}
\def\PO{\operatorname{PO}}
\def \Pr{{\rm Pr}}
\def\PO{\operatorname{PO}}
\def\PSL{\operatorname{PSL}}

\def \R{{\mathbb R}}
\def \Re{{\rm Re}}

\def \sm{\smallsetminus}
\def \SL{{\rm SL}}
\def \SO{{\rm SO}}
\def \st{{\rm st}}
\def\supp{\operatorname{supp}}

\def \tr{{\rm tr\,}}

\def \ul{\underline}

\def\what{\widehat}
\def \Z{{\mathbb Z}}
\def \({\left(}
\def \){\right)}

\newcommand{\sfrac}[2]
{{\textstyle \frac #1 #2}}   

\newcommand{\stack}[2]
{\genfrac{}{}{0pt}{1}{#1}{#2}}

\newcommand{\norm}[1]
{|\hspace{-1pt}| #1 |\hspace{-1pt} |}

\renewcommand{\sp}[1]
{\left\langle #1\right\rangle}

\newcommand{\smat}[4]
{\(\begin{smallmatrix}#1 & #2 \\ #3 & #4\end{smallmatrix}\)}

\newcommand{\esmat}[4]
{\left[\begin{smallmatrix}#1 & #2 \\ #3 & #4\end{smallmatrix}\right]}

\newcommand{\emat}[4]
{\left[\begin{matrix}#1 & #2 \\ #3 & #4\end{matrix}\right]}

\newtheorem{theorem}{Theorem}[section]

\newtheorem{lemma}[theorem]{Lemma}

\newtheorem{proposition}[theorem]{Proposition}

\begin{document}

\pagestyle{myheadings} \markright{FOURIER EXPANSION}

\title{Fourier expansion along geodesics on Riemann surfaces\\ \ \\ \small
Central European Journal of Mathematics 12, No 4, 559-573 (2014)}
\author{Anton Deitmar\\ \ \\
Institute of Mathematics\\ University of Tuebingen\\ Auf der Morgenstelle 10\\ 72076 Tuebingen, GERMANY.
}
\date{}
\maketitle

{\bf Abstract.}
For an eigenfunction of the Laplacian on a hyperbolic Riemann surface, the coefficients of the Fourier expansion are described as intertwining functionals. All intertwiners are classified. A refined growth estimate for the coefficients is given and a summation formula is proved.

\tableofcontents

\section*{Introduction}
For an automorphic function, the invariance under parabolic elements is used to give the standard Fourier expansion, the coefficients of which define the L-function of the form.
In this paper, we instead consider the Fourier-expansion along a hyperbolic element.
In other terms, let $Y$ be a hyperbolic Riemann surface and let $c$ be a closed geodesic in $Y$. 
We are interested in the Fourier coefficients
$$
c_k(f)=\int_0^1f(c(l_ct)e^{2\pi ikt}\,dt
$$
of a smooth function $f\in C^\infty(Y)$. 
Here $l_c$ is the length of the geodesic $c$. 
Under the assumption that $f$ be an eigenfunction of the Laplace operator on $Y$ with eigenvalue $\al$, one can relate $c_k$ to an intertwining integral $I_k^\al(f)$, which depends on $\al$ and $f$, but not on $c$. 
There is an automorphic coefficient $a_k\in\C$, such that
$$
c_k = a_k I_k^\al.
$$
The present paper contains three main results:
\begin{itemize}
\item in Theorem \ref{thm2.3} one finds a classification of all intertwining functionals on the dual of the group $\PGL_2(\R)$.
\item In Theorem \ref{thm2.5} there is given the growth estimate
$$
a_k= O(|k|^\frac12)
$$
as $|k|\to\infty$. The proof uses the technique of analytic continuation developed by Bernstein and Reznikov in \cite{BRan}.
\item In Theorem \ref{thm4.2} finally,  a summation formula is proved, which involves the coefficients $a_k$ and the spectral decomposition in the compact case.
The proof relies on the uniqueness of invariant trilinear forms as in \cite{BR}.
The sum formula is of the form
$$
\sum_{k\in\Z}|a_k|^2\hat\al(k)=\sum_jc_j\int_{\R^2}W_j(t,x)\al(\hat t_x)\,dt\,dx,
$$
where $\al$ is a test function, the decomposition of the $G$-representation on $L^2(\Ga\bs G)$ is $\bigoplus_j\pi_j$ and the constants $c_j$ and the explicit functions $W_j$ depend on $\pi_j$.
Finally $\hat t_x=\frac12\log\left|\frac{(e^{2t}+x)(x-1)}{(e^{2t}+x-1)x}\right|$.
It is hoped that the choice of specific test functions will lead to more precise growth estimates for the $a_k$.
\end{itemize}

We  explain the construction of the factors $a_k$ in a bit more detail. 
Let $X$ be the universal covering of $Y$ and $\Ga$ its fundamental group. 
Then $\Ga$ acts on $X$ by isometries and $Y$ is the quotient $\Ga\bs X$. 
So $\Ga$ injects into the isometry group $G$ of $X$, which acts transitively on $X$, i.e.,  $X\cong  G/K$ for a
maximal compact subgroup $K$. 
Let $(\pi,V_\pi)$ be an irreducible unitary representation of the group $G$ and let 
$\eta:V_\pi \to L^2(\Ga\bs G)$ be an isometric linear $G$-map. 
Let $P_K : L^2(\Ga\bs G)\to L^2(\Ga\bs G)^K = L^2(\Ga\bs G/K) = L^2(Y)$ denote the orthogonal projection onto the subspace of $K$-invariants. 
Demanding that $f \in C^\infty(Y)$ be an eigenfunction of the Laplacian amounts to the same as demanding $f$ to lie in the image of $P_K\circ\eta$ for some $\pi$ and some $\eta$. The functional $I_k^\ga= c_k\circ P_K\circ\eta$ on $V_\pi$ then has an intertwining property with respect to a split torus $A$ inside $G$. 
By a uniqueness result, proven in Section \ref{sect2}, this implies that $I_k^\ga$ is a multiple of a standard intertwiner $I_{\pi,k}^\st$ on $V_\pi$, which we named $I_k^\al$ above. 
So we get the existence of a factor $a_k\in\C$ with $I_k^\ga = a_kI_{\pi,k}^\st = a_kI_k^\al$ as above.

In Section \ref{sect1} we describe the setting in greater precision.
In Section \ref{sect2} we classify the intertwining functionals that show up in the context and define the standard intertwiners that give rise to the factors $a_k$ above.
We also show the growth estimate of the factors $a_k$.
In Section \ref{sect3} we show how the Fourier expansion along a geodesic expands to an expansion on the whole space and in Section \ref{sect4} we show how to derive the summation formula from the uniqueness of triple products.

\section{Generalized period integrals}\label{sect1}
In this paper we use the group $\GL_2(\R)$, the elements of which we write as matrices $\smat abcd$ and the group $G=\PGL_2(\R)=\GL_2(\R)/\R^\times$, the elements of which we write in the form $\esmat abcd$.
Here we will usually arrange the determinant of $\smat abcd$ to be $1$ or $-1$.
The connected component of $G$ is $G^0=\PSL_2(\R)=\SL_2(\R)/(\pm 1)$.
The group $G^0$ acts on the upper half plane $\H$ in $\C$ by linear fractionals and this action extends to an action of $G$ in a way that $G$ is identified with the group of all hyperbolic isometries on $\H$.
The stabilizer in $G$ of the point $i\in\H$ is the maximal compact subgroup $K=\PO(2)={\rm O}(2)/\pm 1$ of $G$.
So $\H$ is identified with $G/K$.

Let $\Ga$ be a discrete subgroup of the group $G$.
Later we will assume $\Ga$ to be of finite covolume.
For simplicity, we will assume $\Ga$ to be torsion-free and that $\Ga\subset G^0$.
This implies that $\Ga$ is the fundamental group of $Y=\Ga\bs\H$ and the latter is a Riemann surface equipped with the hyperbolic metric.

For a closed geodesic $c$ in $Y$ let $l(c)$ denote its length.
The period integral $I_c(f)=\int_0^{l(c)}f(c(t))\, dt$ 
is the zeroth coefficient of the Fourier-expansion of 
the function $t\mapsto f(c(t))$.
Therefore the higher coefficients can be viewed as ``generalised period integrals''.

The real Lie-algebra of $G$ is $\g_\R=sl_2(\R)$, the Lie-algebra of all real $2\times 2$ matrices of trace zero.
For $X,Y\in\g_\R$ let $b(X,Y)=\frac12\tr(XY)$.
Then $b$ is an invariant 
symmetric bilinear form.
Let $\k_\R\subset \g_\R$ be the Lie algebra of $K$.
Then $b$ is negative definite on $\k_\R$ and positive 
definite on its orthogonal complement $\p_\R$.
Let $\a_\R=\R\smat 1\ \ {-1}$, and 
let $A=\exp(\a_\R)$ be the corresponding subgroup of $G$.
Let $\a_\R^+=\R_{>0}\smat 1\ \ {-1}$ be the positive cone and let $A^+=\exp(\a_\R^+)$.
Then $A$ is closed and non-compact and its centralizer in $G$ is the group $AM$ of all diagonal matrices in $G$.
Here $M$ is the two element group generated by $\esmat{-1}\ \ 1$.
Let $\a$ be the complexification of $\a_\R$ and $\a^*$ the dual space of $\a$.
Then $\a^*$ can be identified with the set of continuous homomorphisms from $A$ to $\C^\times$.
For $\la\in\a^*$ we write $a\mapsto a^\la=e^{\la(\log a)}$ for the corresponding homomorphism.
It is known that $G^0\to G/K$ can be identified with the 
sphere-bundle of $\H=G/K$ in a way that the geodesic flow 
is given by
$$
\phi_t(g)= g\exp(tH_1),\qquad g\in G^0,
$$
where $H_1=\esmat 1\ \ {-1}$.
The sphere bundle $SY$ equals $\Ga\bs G^0$.

A closed geodesic $c$ in $Y$ gives rise to a conjugacy 
class $[\ga]$ in $\Ga$ of elements which ``close'' $c$.
Any such $\ga\in \Ga$ is \emph{hyperbolic} in the sense 
that it is conjugate in $G$ to an element of the form $a_{t_0}=\diag(e^{t_0},e^{-t_0}) \in A\sm\{ 1\}$.
We insist that $a_{t_0}\in A^+$, i.e., $t_0>0$ to make it unique.
We now assume that $\ga$ be primitive, i.e., $\ga$ is no power $\tau^n$ for any $\tau\in\Ga$ and $n\ge 2$.
This is equivalent to the geodesic $c$ being primitive, i.e., $c$ is no power of any shorter geodesic.
The characters of the compact abelian group $A/\sp{a_\ga}$ are given by those $\mu\in\a^*$ with $a_\ga^\mu=1$, i.e., $\mu(\log a_\ga)\in 2\pi i\Z$.
Let $\mu_{\ga}$ be the unique element of $\a^*$ with $\mu_\ga(\log a_\ga)=2\pi i$.
Then $\widehat{A/\sp{a_\ga}}=\Z \mu_\ga$.
Later we will use the notation
$$
\tilde\mu=\frac1{2\pi i}\mu.
$$

Fix an element $\sigma=\sigma_\ga\in G^0$ with $\ga = \sigma^{-1} a_{t_0}\sigma$.
Let $f\in C^\infty(\Ga\bs G)=C^\infty(\Ga\bs G)$ and set
$$
f^\sigma (x)=f(\sigma x).
$$
Then the map $t\mapsto f^\sigma(a_t x)$ with $a_t=\diag(e^t,e^{-t})$ is periodic of period $t_0$ and thus has a Fourier-expansion
$$
f^\sigma(a_tx)=\sum_{k\in\Z}e^{2\pi i t/t_0}\frac1{t_0}\int_0^{t_0}f^\sigma(a_tx)e^{-2\pi i t/t_0}\,dt.
$$
For $k\in\Z$ let
\begin{align*}
I_k^\ga : C^\infty(\Ga\bs G) &\to  \C\\
f &\mapsto  \frac 1{t_0}\int_0^{t_0}f^\sigma(a_t)e^{-2\pi i kt/t_0}\,dt.
\end{align*}
Note that $I_k^\ga$ depends on the choice of $\sigma$. 
Geometrically, this corresponds to choosing a base-point on the closed orbit $c$.
This dependence is not severe, as $\sigma$ is determined up to 
multiplication from the right by elements of $A
$.
If we replace $\sigma$ by $\sigma a_0$, then 
$I_k^\ga$ is replaced by $a_0^{k\mu_\ga}I_k^\ga$.
So in particular, the absolute value $|I_k^\ga|$ is uniquely determined by $k$ and $\ga$.
Further, if $\ga$ is replaced by a $\Ga$-conjugate, say $\ga'=\tau\tau^{-1}$ then one can choose $\sigma_{\ga'}$ to be equal to $\tau \sigma_\ga$ and with this choice one gets $I_k^{\ga'}=I_k^\ga$.

The form $b$ determines a Haar measure $dg$ on $G$.
Let $R$ denote the unitary $G$-representation on $L^2(\Ga\bs G)$ given by right translations.
We are particularly interested in the subspace $L^2_\cusp(\Ga\bs G)$ of cusp forms.
This representation space decomposes discretely,
$$
L^2_\cusp(\Ga\bs G)\ \cong\ \bigoplus_{\pi\in\hat G} N_\Ga(\pi)\pi,
$$
where the sum runs over the unitary dual $\hat G$ of $G$ and the multiplicities $N_\Ga(\pi)$ are finite.
Here and later we understand the direct sum to be a completed direct sum in the appropriate topology.
We define $C^\infty_\cusp(\Ga\bs G)$ to be the intersection of $C^\infty(\Ga\bs G)$ with the space of cusp forms.
Then it turns out that $C_\cusp^\infty(\Ga\bs G)$ is the set of smooth vectors in the representation space $L_\cusp^2(\Ga\bs G)$, i.e., 
$$
C^\infty_\cusp(\Ga\bs G)= L^2_\cusp(\Ga\bs G)^\infty= \bigoplus_{\pi\in\hat G} N_\Ga(\pi)\pi^\infty,
$$
where $\pi^\infty$ is the representation on the Fr\'echet space of smooth vectors.
 The linear functional $I_k^\ga$ satisfies
$$
I_k^\ga(R(a)\ph)= a^{k\mu_\ga} I_k^\ga(\ph)
$$
for every $a\in A$.
This means that $I_k^\ga$ is an intertwining functional.

\section{Intertwining functionals}\label{sect2}
We first shall give 
a description of the admissible and unitary duals of the group $G$.
For any topological group $G$, let $\what G$ denote the \e{unitary dual} of $G$, that is, the set of unitary equivalence classes of irreducible unitary representations of $G$.

Next let $G$ denote a semisimple Lie group with finite center and finitely many connected components.
Then $G$ has a maximal compact subgroup $K$ which is unique up to conjugation.
A representation $(\pi,V_\pi)$ of $G$ is called \e{admissible}, if for every $\tau\in\what K$ the isotype $V_\pi(\tau)$ is finite-dimensional.
In that case the space $V_{\pi,K}$ of $K$-finite vectors in $V_\pi$ forms a $(\g,K)$-module, where $\g$ is the complexified Lie algebra of $G$.
Two admissible representations are called \e{infinitesimally equivalent} if their $(\g,K)$-modules of $K$-finite vectors are isomorphic.
The \e{admissible dual} $\what G_{\adm}$ of $G$ is the set of infinitesimal equivalence classes of irreducible admissible representations of $G$.

Due to results of Harish-Chandra, every irreducible unitary representation of $G$ is admissible and two unitary admissible representations are unitarily equivalent if and only if they are infinitesimally equivalent.
Thus the unitary dual $\what G$ can be considered a subset of the admissible dual $\what G_\adm$.

Now consider $G=\PGL_2(\R)$.
There is a canonical character
$$
\chi: G\to\{\pm 1\};\qquad g\mapsto {\rm sign}(\det(g)),
$$
taking the values $\pm 1$ and having the connected component $G^0$ for kernel.
For a representation $\pi$ of $G$ we define the \e{$\chi$-twist} $\chi\pi$ of $\pi$, also written $\chi\otimes\pi$ as the representation with the same space $V_\pi$ as $\pi$ but defined as
$$
\chi\pi(x)=\chi(x)\pi(x).
$$

Let $P$ denote the parabolic subgroup of $G$ consisting of all upper triangular matrices.
Then $P=MAN$, where $N$ is the group of all upper triangular matrices with ones on the diagonal.
For $\la\in\C$ and $a=\esmat{e^t}\ \ {e^{-t}}\in A$ we write
$$
a^\la=e^{\la t}.
$$

Let $\pi_{\la}$ be the corresponding principal series representation, which we normalize to live on the space of functions $\ph:G\to \C$ satisfying $\ph(manx)=a^{\la+1} \ph(x)$.
Let $V_{\la}$ be the space of $\pi_{\la}$ and let $V_{\la}^\infty$ be the space of smooth vectors in it.
These can be viewed as smooth sections of the line bundle $E_{\la}$ over $P\bs G$ given by the $P$-representation $(man)\mapsto a^{\la+1}$.
Note that restriction of functions to $K$ identifies $V_{\la}^\infty$ with the space $C^\infty(M\bs K)$, so in particular, for $\ph\in V_\la^\infty$, the function $\ph|_K$ is independent of $\la$.

\begin{proposition}
Let $G=\PGL_2(\R)$.
The admissible dual $\what G_\adm$ of $G$ consists of
\begin{enumerate}[\rm (a)]
\item $\pi_\la$, $\la\in\C$, $\la\notin 1+ 2\Z$,
\item $\CD_{2n}$ for $n=1,2,3,\dots$ the standard discrete series representations,
\item $\delta_m$ for $m=0,2,4,\dots$ where $\delta_m$ is the $(m+1)$-dimensional representation on the space of all homogeneous polynomials $p(X,Y)$ of degree $m$,  
\end{enumerate}
together with their $\chi$-twists.
The only isomorphisms occurring are $\pi_{it}\cong\pi_{-it}$ and $\chi\pi_{it}\cong\chi\pi_{-it}$ for $t\in\R$.

The unitary dual $\what G$ consists of
\begin{enumerate}[\rm (a)]
\item $\pi_\la$ for $\la\in i[0,\infty)\cup(0,1)$,
\item $\CD_{2n}$ for $n=1,2,3,\dots$ the discrete series representations,
\item $\delta_0$, 
\end{enumerate}
together with their $\chi$-twists.

For $n=1,2,3,\dots$ we have the exact sequence of representations
$$
0\to \CD_{2n}\to\pi_{2n-1}\to\delta_{2n-2}\to 0,
$$
and
$$
0\to\delta_{2n-2}\to\pi_{1-2n}\to\CD_{2n}\to 0.
$$
\end{proposition}

\begin{proof} This can be deduced from the description of the unitary dual of $\SL_2(\R)$ given, for example, in \cite{Knapp}.
\end{proof}

Let $T=\esmat{-1}\ \ 1$ be the non-trivial element of $M$.
As $T^2=1$, for every representation $(\pi,V_\pi)$, the space $V_\pi$ splits as a direct sum $V_\pi=V_\pi^+\oplus V_\pi^-$, where $V_\pi^\pm$ is the $\pm 1$-eigenspace of $\pi(T)$.
Twisting by $\chi$ interchanges the roles of $V_\pi^+$ and $V_\pi^-$.

Let $\mu\in\a^*$ and let $(\pi,V_\pi)$ be a representation of $G$.
A continuous linear functional $l:V_\pi^\infty\to\C$ is called a \e{$\mu$-intertwiner}, if
$$
l(\pi(a)v)= a^{\mu} l(v)
$$
holds for every $a\in A$ and every $v\in V_\pi$.
Let $V_\pi^\infty(\mu)$ be the space of all $\mu$-intertwiners.
Note that $V_\pi^\infty(\mu)=V_{\chi\pi}^\infty(\mu)$, where we consider the $\chi$-twist $\chi\pi$ as a representation with the same representation space as $\pi$.

Let $w_0=\esmat\ {-1}1\ $ and $n_0=\esmat 11\ 1$.
Then $w_0$ is a representative of the non-trivial element of the Weyl-group $W(G,A)$.

The base space $P\bs G$ of the bundle $E_{\la}$ consists of three orbits under the group $AM$,
namely the open orbit $[w_0n_0]$, and the two closed orbits $[1]$, $[w_0]$ which are indeed points.

We now define a standard intertwiner on the representation $\pi_\la$ for $\Re(\la)>-1$.
Let
$$
I_{\la,\mu}^{\st}(\ph)=I^{\st}_{\pi_\la,\mu}(\ph)= \int_A \ph(w_0n_0 a)a^{- \mu}\, da,\qquad \ph\in V_\la^\infty.
$$
If $\supp \ph\subset [w_0n_0]$, then the integral $I_{\la,\mu}^{\st}$ is extended over a compact set, hence it exists.

\begin{lemma}
Let $\mu\in\C$.
If $\Re(\la)>-1-\Re(\mu)$, then the integral $I_{\la,\mu}^{\st}(\ph)$ exists for every $\ph\in V_\la^\infty$ and defines a $\mu$-intertwiner.
The map $\la\mapsto I_{\la,\mu}^{\st}$ extends to a meromorphic operator-valued function with poles exactly at 
$$
\la= -\mu-1-4k,\qquad k\in\N_0
$$
and
$$
\la=\mu-3-4k,\qquad k\in\N_0.
$$
Outside the poles, $I_{\la,\mu}^{\st}$ spans the one dimensional space of intertwiners supported on the open orbit.
At the poles this space is zero.
\end{lemma}

\begin{proof}
Let $\ul a:G\to A$, $\ul n: G\to N$, $\ul k:G\to K$ be the smooth maps defined by the Iwasawa decomposition $g=ank=\ul a(g)\ul n(g)\ul k(g)$ for $g\in G$.
For $\ph\in V_\la^\infty$ one has
\begin{align*}
I_{\la,\mu}^{\st}(\ph) &= \int_A\ul a(w_0n_0a)^{\la+1} \ph(\ul k(w_0n_0a)) a^{-\mu}\,da\\
&=\int_A \ul a(w_0n_0^a)^{\la+1}a^{-\la-1-\mu}\ph(\ul k(w_0n_0^a))\,da,
\end{align*}
where $n_0^a=a^{-1}n_0a$.
Noting that if $a=\esmat{e^t}\ \ {e^{-t}}$, then we have $n_0^a=\esmat 1{e^{-2t}}\ 1$, we get $w_0n_0^a=\esmat 0{-1}1{e^{-2t}}$ and thus $\ul k(w_0n_0^a)=\frac1{\sqrt{e^{-4t}+1}}\esmat{e^{-2t}}{-1}1{e^{-2t}}$.
We now define special test functions.
For $k,l\in\N_0$ let $\ph_{k,l}:\SO(2)\to\C$ be defined by
$$
\ph_{k,l}\smat d{-c}cd=c^k d^l.
$$
If $k+l$ is even, $\ph_{k,l}$ defines an element of $V_\la^\infty$.
The theory of Taylor-series tells us that every $\ph\in V_\la^\infty$ can be written as
$$
\ph=\sum_{{0\le k,l\le N}}c_{k,l}\ph_{k,l}+R_N(\ph),
$$
where $c_{k,l}\in\C$ and the function $R_N(\ph)\in V_\la^\infty\cong C^\infty(M\bs K)$ vanishes to order $N$ at $1$ and $w_0$.
The integral $I_{\la,\mu}^{\st}(\ph)$, as written above, makes sense for $\ph=\ph_{k,l}$ also in the case when $k+l$ is not even.
We use this fact for convenience.
For $k,l\in\N_0$ we compute
\begin{align*}
I_{\la,\mu}^{\st}(\ph_{k,l}) &= \int_A \ul a(w_0n_0^a)^{\la+1}a^{-\la-1-\mu}\ph_{k,l}(\ul k(w_0n_0^a))\,da\\
&= \int_\R \ul a\smat 0{-1}1{e^{-2t}}^{\la+1}e^{-(\la+1+\mu) t}\ph_{k,l}\(\ul k\smat 0{-1}1{e^{-2t}}\)\, dt\\
&=\int_\R (e^{-4t}+1)^{-\frac{\la+1}2}
e^{-(\la+1+\mu) t}
\ph_{k,l}\(\frac1{\sqrt{e^{-4t}+1}}\smat{e^{-2t}}{-1}1{e^{-2t}}\)\,dt\\
&=\int_\R (e^{-4t}+1)^{-\frac{\la+k+l+1}2}
e^{-(\la+2l+1+\mu) t}\,dt\\
&=\frac14\int_\R (e^{-t}+1)^{-\frac{\la+k+l+1}2}
e^{-\frac{\la+2l+1+\mu}4 t}\,dt\\
&= \frac14B\(\frac{\la+2l+1+\mu}4,\frac{\la+2k+1-\mu}4\),
\end{align*}
where $B(x,y)=\frac{\Ga(x)\Ga(y)}{\Ga(x+y)}$ is Euler's Beta-function.
We conclude
$$
I_{\la,\mu}^{\st}(\ph_{k,l})=\frac{\Ga\(\frac{\la+2l+1+\mu}4\)\Ga\(\frac{\la+2k+1-\mu}4\)}{\Ga\(\frac{\la+k+l+1}2\)}.
$$
Note that
$$
I_{\la,\mu}^{\st}(\ph_{k+1,l+1})=I_{\la+2,\mu}^{\st}(\ph_{k,l}).
$$
The space $\ph_{N,N}C^\infty(M\bs K)$ is the space of all $\ph\in V_\la^\infty$ which vanish to order $\ge N$ at $1$ and $w_0$.
For any $\ph$ in this space, the integral $I_{\la,\mu}^{\st}(\ph)$ converges if $\Re(\la)>-N-1-\Re(\mu)$.
Therefore we get analytic continuation of the map $\la\mapsto I_{\la,\mu}^{\st}$ as claimed.
The lemma follows.
\end{proof}

We next consider intertwiners which are supported on the closed orbits $[1]$ and $[w_0]$.
Let $S_{[1],0}: V_\la^\infty\to\C$ denote the distribution
$$
S_{[1],0}(\ph)= \ph(1).
$$
Then $S_{[1],0}\circ R(a)=a^{\la+1}S_{[1],0}$, so $S_{[1],0}$ is an $\mu$-intertwiner for $\mu={\la+1}$.

We next consider higher derivatives of this distribution.
For $X\in\g$, the Lie algebra of $G$, and $f\in C^\infty(G)$, we let
$$
R_Xf(y)=\left.\frac d{dt}\right|_{t=0}f(y\exp(tX)).
$$

Let $\bar N=\theta(N)$ and let $\bar\n_\R$ be its Lie algebra.
Then the tangent space of $P\bs G$ at the unit is isomorphic to $\bar n_\R$.
Let $X=\smat 0100 $ and $\bar X=\smat 0010$. Then $\n_\R=\R X$ and $\bar\n_\R=\R \bar X$.
For $k\in\N$ set
$$
S_{[1],k}(\ph)= R_{\bar X}^k \ph(1).
$$
Then $S_{[1],k}\circ R(a)=a^{\la+1+{2k}}S_{[1],k}$.
Since these span the space of all distributions supported at $1$ we see that we get a non-zero $\mu$-intertwiner supported on $1$ if and only if $$
\mu=\la+1+{2k}
$$ 
for some $k\in\N_0$.
If this condition is satisfied, then the space of intertwiners supported on $[1]$ is one dimensional.

We turn to the other closed orbit $[w_0]$.
In this case we define
$$
S_{[w_0],k}(\ph)= R_X^k \ph(w_0).
$$
Then $S_{[w_0],k}\circ R(a)=a^{-\la-1-{2k}}S_{[w_0],k}$ and we conclude that there exists a non-zero intertwiner supported on $[w_0]$ if and only if
$$
\mu=-\la-1-2k,\qquad k\in\N_0,
$$
in which case the space of intertwiners supported on $[w_0]$ is one dimensional.

\begin{theorem}\label{thm2.3}
\begin{enumerate}[\rm (a)]
\item For $\la,\mu\in\C$ with $\mu\ne 0$ we have
$$
\dim V_\la^\infty(\mu)=1.
$$
In this case, if $\la\notin \pm\mu-1-2\N_0$, then $V_\la^\infty(\mu)$ is spanned by $I_{\la,\mu}^{\st}$.
If $\la=\mu-1-2k$ with $k\in\N_0$, then $V_\la^\infty(\mu)$ is spanned by $S_{[1],k}$ and if
$\la=-\mu-1-2k$ with $k\in\N_0$, then $V_\la^\infty(\mu)$ is spanned by $S_{[w_0],k}$.

The same holds for the $\chi$-twist.
 
\item For $\la\in\C$ and $\mu=0$ we have
$$
\dim V_\la^\infty(0)=\begin{cases}
2& \la\in-1-2\N_0,\\
1&\text{otherwise}.
\end{cases}
$$
If $\la\notin -1-2\N_0$, then $V_\la^\infty(0)$ is spanned by $I_{\la,0}^{\st}$.
If $\la=-1-2k$ with $k\in\N_0$, then $V_\la^\infty(0)$ is spanned by $S_{[w_0],k},S_{[1],k}$.

The same holds for the $\chi$-twist.

\item
For the finite-dimensional representations we have
$\dim\delta_{2n}(\mu)= 0$ if $\mu\ne -2n,-2n+2,\dots,2n$ and
$$
\dim\delta_{2n}(\mu)=1\quad\text{if }\mu\in\{-2n,-2n+2,\dots,2n\}.
$$

\item
Let $\mu\in\C$, $n\in\N$. Then we have an  exact sequence
$$
0\to\CD_{2n}^\infty(\mu)\to V_{1-2n}^\infty(\mu)\to\delta_{2n-2}(\mu)\to 0.
$$
If $\mu\notin \{-2n+2,\dots,2n-2\}$, then $\delta_{2n-2}(\mu)=0$ and therefore $\CD_{2n}^\infty(\mu)\cong V_{1-2n}^\infty(\mu)$ is one dimensional.

If $\mu\in \{-2n+2,\dots,2n-2\}$ but $\mu\ne 0$, then $\CD_{2n}^\infty(\mu)=0$ and finally $\CD_{2n}^\infty(0)$ is one-dimensional and is spanned by $S_{[1],n-1}-S_{[w_0],n-1}$.
\end{enumerate}\end{theorem}

\begin{proof}
(a) and (b) are clear by the above.
For (c) recall that $\delta_{2n}$ has the basis $e_{2j-2n}= X^jY^{2n-j}$ for $j=0,\dots 2n$ and the group $A$ acts by $\delta_{2n}(a)e_{2j-2n}= a^{2j-2n}e_{2j-2n}$.
This proves (c).

For (d) we consider the exact sequence
$$
0\to\delta_{2n-2}\to\pi_{1-2n}\to\CD_{2n}\to 0,
$$
which induces the exact sequence of intertwiners
$$
0\to\CD_{2n}^\infty(\mu)\to V_{1-2n}^\infty(\mu)\to\delta_{2n-2}(\mu).
$$
This proves the first assertion, i.e., the case $\delta_{2n-2}(\mu)=0$.
If $\mu\in \{-2n+2,\dots,2n-2\}$, which means $\delta_{2n-2}(\mu)\ne 0$, ten we have to show that the map $V_{1-2n}^\infty(\mu)\to\delta_{2n-2}(\mu)$ is non-zero, for it is automatically onto then, as the target space is one-dimensional.
The above exact sequence dualizes to the exact sequence
$$
0\to \CD_{2n}\to \pi_{2n-1}\to \delta_{2n-2}\to 0,
$$
which yields an exact sequence
$$
0\to \delta_{2n-2}(\mu)\to V_{2n-1}^\infty(\mu)\to \CD_{2n}^\infty(\mu).
$$
So the arrow $\delta_{2n-2}(\mu)\to V^\infty_{2n-1}(\mu)$ is non-zero, hence its dual $V_{1-2n}^\infty(\mu)\to\delta_{2n-2}(\mu)$ likewise.

Finally, for $\mu=0$ we show that the kernel of the restriction map 
$V_{1-2n}^\infty(0)\to\delta_{2n-2}(0)$ is spanned by $S_{[w_0],n-1}-S_{[1],n-1}$.
Recall that $R_{X}$ and $R_{\bar X}$ are the weight-change operators in highest weight theory, so $R_X^{n-1}e_0$ is a multiple of $e_{2n-2}$ and $R_{\bar X}^{n-1}e_0$ is a multiple of $e_{2-2n}$.
Let $\theta(x)=x^{-t}$ be the transpose followed by the inversion, i.e., $\theta$ is the Cartan-involution on $G$ with fixed point set $K$.
We have $\theta(X)=\bar X$.
It is readily verified that the map $\Psi$ with $\Psi(f)(x)=f(w_0\theta(x))$ is an isomorphism between the representations $\pi_\la$ and $\pi_\la\circ\theta$.
Thus it follows that $S_{[w_0],n-1}(\Psi(\ph))=S_{[1],n-1}(\ph)$.
As $\Psi(\pi(a)\ph)=\pi(a^{-1})\Psi(\ph)$ and $\Psi^2=\Id$ it follows that $\Psi(e_0)=\pm e_0$.
We want to show $\Psi(e_0)=e_0$.
For this we let $\smat d{-c}cd\in \SO(2)$ with $c,d>0$.
The $A$-invariance of $e_0$ shows that for $a>0$ we have
$$
e_0\(\frac1{\sqrt{a^2c^2+d^2/a^2}}\smat {ad}{-c/a}{ac}{d/a}\)=(c^2a^2+d^2/a^2)^{\frac{2-n}2}e_0\smat d{-c}cd.
$$
So the $A$-orbit is mapped to positive multiples of $e_0\smat d{-c}cd$.
By the $M$-invariance we get on the other hand,
$$
e_0\(\smat\ {-1}1\ \smat d{-c}cd\)=e_0\smat{-c}{-d}d{-c}=e_0\smat c{-d}dc.
$$
As $\smat c{-d}dc$ lies in the same $A$-orbit as $\smat d{-c}cd$, the claim follows.
We therefore conclude that $S_{[w_0],n-1}(e_0)=S_{[1],n-1}(e_0)$.
\end{proof}

A sequence $(c_j)_{j\in\N}$ of complex numbers is said to be \emph{of moderate growth}, if there exist $N\in\N$ such that
$$
|c_j|= O(j^N),
$$
as $j\to\infty$.
The sequence is called \emph{rapidly decreasing}, if for every $N\in\N$ one has
$$
|c_j|= O(j^{-N})
$$
as $j\to\infty$.
The product of two moderately growing sequences is moderately growing and the product of a moderately growing sequence and a rapidly decreasing sequence is rapidly decreasing.

\begin{proposition}
Let $\Re(\la)>-1$ and $f\in V_\la^\infty$ as well as $\mu\in i\R\sm\{0\}$ be given.
Then the sequence $(I_{\la,k\mu}^{\st}(f))_{k\in\Z}$ is rapidly decreasing.
\end{proposition}

\begin{proof}
Recall $H_1=\smat 1\ \ {-1}$.
Using integration by parts, we compute for $N\in\N$,
\begin{align*}
I_{\la,k\mu}^{\st}(\ph) 
&= \int_A \ph(w_0n_0a)a^{-k\mu}\,da\\
&=\int_\R \ph(w_0n_0\exp(tH_1)e^{-k\mu t}\,dt\\
&=\(\frac{1}{k\mu}\)^N\int_\R R_{H_1}^N \ph(w_0n_0a)a^{-k\mu}\,da.
\end{align*}
As $\Re(\mu)=0$ we have $|a^{-k\mu}|=1$ and the claim follows.
\end{proof}

By an \emph{automorphic representation} $(\pi,V_\pi,\eta)$ we mean an irreducible unitary representation $(\pi,V_\pi)$ of $G$ together with an isometric $G$-equivariant linear map $\eta\colon V_\pi\to L^2(\Ga\bs G)$.
Then $\eta$ maps the space $V_\pi^\infty$ of smooth vectors into $C^\infty(\Ga\bs G)$.
Let $k\in\Z$.
The map $I_k^\ga\circ \eta$ is an intertwiner for $k\mu_\ga$.
So, for instance, let $\pi=\pi_\la\in\what G$, then the space of intertwiners on $V_\la$ is spanned by $I_{\la,k\mu_\ga}^{\st}$.
Therefore there exists $a_k^{\eta,\ga}\in\C$ such that
$$
I_k^\ga\circ \eta=a_k^{\eta,\ga} I_{\la,k\mu_\ga}^{\st}.
$$
However, if $\pi=\CD_{2n}$ is a discrete series representation and $k=0$, then there are $b_0^{\eta,\ga},c_0^{\eta,\ga}\in\C$ such that
$$
I_0^{\ga}\circ\eta= b_0^{\eta,\ga} S_{[1],n-1}+c_0^{\eta,\ga}S_{[w_0],n-1},
$$
where $c_0^{\eta,\ga}=-b_0^{\eta,\ga}$ if $n$ is even and analogously for the $\chi$-twist.
Finally, If $\pi$ is the trivial representation, 
we consider $\pi$ as a subrepresentation of $\pi_{-1}$, so this case does not need extra treatment.

\begin{theorem}\label{thm2.5}
\begin{enumerate}[\rm (a)]
\item There exists a constant $C_\ga>0$, depending on $\ga$, such that
$$
|a_k^{\eta,\ga}|\le C_\ga \(1+|\la|^{\frac14}\)\(1+|\la^2-k^2\mu^2|^{\frac14}\)e^{\frac\pi 4|k\mu|}
$$
holds for every cuspidal automorphic representation $\eta: V_\la\to L^2_\cusp(\Ga\bs G)$ and every $k\in\Z$.

\item For a fixed automorphic representation $\eta: V_\la\to L_\cusp^2(\Ga\bs G)$ there exists a constant $C_{\eta,\ga}>0$ such that
$$
|a_k^{\eta,\ga}|\le C_{\eta,\ga} (1+|k|^\frac12)
$$
holds for every $k\in\Z$.

\item There exists a constant $D_\ga>0$, depending on $\ga$, such that
$$
|b_0^{\eta,\ga}|, |c_0^{\eta,\ga}|\le D_\ga n^{5-n}
$$
holds for every cuspidal automorphic representation $\eta:\CD_{2n}\to L^2_\cusp(\Ga\bs G)$.
\end{enumerate}
\end{theorem}

\begin{proof}
(a)
Let $\ph=\ph_\la\in V_{\la}$ be the unique $K$-invariant function which on $K$ takes the value $1$ and let $f=\eta(\ph)$.
Then
\begin{align*}
I_k^\ga(f)&= a_k^{\eta,\ga}I_{\la}^{k\mu}(\ph)= a_k^{\eta,\ga}I_{\la}^{k\mu}(f_{0,0})=a_k^{\eta,\ga}
\frac{\Ga\(\frac{\la+1+k\mu}4\)\Ga\(\frac{\la+1-k\mu}4\)}{\Ga\(\frac{\la+1}2\)}
\end{align*}
If $\Ga$ is cocompact, by \cite{SeegerSogge} the sup norm of $f$ satisfies
$$
\norm{f}_\infty= O\(|\la|^{1/4}\).
$$
If $\Ga$ is not cocompact, one finds in \cite{Iwaniec}, Sec. 13.2, that
$$
|f(z)|=O_z(|\la|^{1/4}),
$$
where the implied constant depends continously on $z$.
Therefore, this estimate holds uniformly on the closed geodesic attached to $\ga$.
According to \cite{GR} 8.328.1, for fixed real $x$ and for $|y|\to\infty$ one has
$$
|\Ga(x+iy)|\sim\sqrt{2\pi}e^{-\frac\pi 2|y|}|y|^{x-\frac12}.
$$
This implies the claim.

(b) We start out as in the above proof, except that we don't use the estimate $|f(z)|\ll_z |\la|^{\frac14}$.
We thus get
$$
|a_k^{\eta,\ga}| \ll |I_k^\ga(f)|(1+|\la^2-k^2\mu^2|^\frac14)e^{\frac\pi 4|k\mu|}.
$$
It suffices to show that for fixed $\ga$ and $\eta$ one has
$$
|I_k^\ga(f)|\ll e^{-\frac\pi 4 |k\mu|}.
$$
We will show this using the technique of analytic continuation of representations from \cite{BRan}.
Let $X\in sl_2(\R)$ with $\ga=\exp(X)$.
After conjugating $\Ga$, we may assume $X=\diag(A,-A)$ for some $A>0$.
Then $\ga=\diag(e^A,e^{-A})$ and $\mu=\frac{2\pi i}A$.
We have
$$
I_k^\ga(f)=\int_0^1f(\exp(tX))e^{-2\pi ikt}\,dt.
$$
It is easy to see that the function $t\mapsto \pi_\la(a_t)\ph)$ with $a_t=\esmat{e^{At}}\ \ {e^{-At}}$ extends to a holomorphic function from $\{ |\Im z|<\frac\pi{4A}\}$ to $V_\la$.
It follows that the function $f(\exp(tX))=\eta(\pi_\la(a_t)\ph)(1)$ extends to a holomorphic function on the set of all $z=x+iy\in\C$ with $|y|<\frac\pi{4A}$.
We get a continuous extension to $|y|\le\frac\pi {4A}$.
For $k\ge 0$ we get by a shift of the contour integral that
\begin{align*}
I_k^\ga(f)&=\int_0^1f(\exp(tX))e^{-2\pi ikt}\,dt\\
&= \int_0^1f(\exp((t-i\frac\pi{4A})X) e^{-2\pi ikt}\,dt\ e^{-\frac{\pi^2}{2A}k}={\rm const}\cdot e^{-\frac\pi 4|k\mu|}.
\end{align*}
For $k<0$ we similarly move the contour to $i\frac\pi{4A}$.

(c)
Let $n\in \N$ and let $\ph_n:K\to\C$ be given by
$$
\ph_n\(\esmat d{-c}cd\)=(ci+d)^{2n}.
$$
Then $\ph_n$ and its complex conjugate span the lowest $K$-type in $\CD_{2n}\cong\pi_{1-2n}/\delta_{2n-2}$.
A computation shows
\begin{align*}
R_{\bar X}\ph_n&=i n(\ph_n+\ph_{n+1}),\\
R_{X}\ph_n&=-i n(\ph_n-\ph_{n+1}).
\end{align*}
We use induction in $k\in\N_0$ to show that
$$
S_{[1],k}(\ph_n)=(-1)^{k+n}S_{[w_0],k}(\ph_n).
$$
For $k=0$ we have
$$
S_{[1],k}(\ph_n)=\ph_n(1)=1=(-1)^n\ph(w_0)=(-1)^{k+n}S_{[w_0],k}(\ph_n).
$$
The step $k\mapsto k+1$ is
\begin{align*}
S_{[1],k+1}(\ph_n)&=R_{\bar X}^{k+1}\ph_n(1)\\
&=in (R_{\bar X}^k(\ph_n+\ph_{n+1})(1)\\
&=in(S_{[1],k}(\ph_n)+S_{[1],k}(\ph_{n+1})\\
&=in((-1)^{k+n}S_{[w_0],k}(\ph_n)+(-1)^{k+n+1}S_{[w_0],k}(\ph_{n+1})\\
&=(-1)^{k+n+1}\(-in(S_{[w_0],k}(\ph_n)+S_{[w_0],k}(\ph_{n+1})\)\\
&=(-1)^{k+n+1}\(-in(R_X^k\ph_n(w_0)+R_X^k\ph_{n+1}(w_0))\)\\
&=(-1)^{k+n+1}\(R_X^{k+1}\ph_n(w_0)\)\\
&=(-1)^{k+n+1}\(S_{[w_0],k+1}(\ph_n)\).
\end{align*}
Next we show that
$$
|S_{[1],{n-1}}(\ph_n)|\ge n^{n-1} \quad\text{and}\quad |S_{[1],{n-1}}(\ph_{n+1})|\ge n^{n-1}
$$
This follows from $\ph_m(1)=1$ and the fact that 
$$
R_{\bar X}^k\ph_n=i^k(n^k\ph_n+(*)),
$$
where $(*)$ denotes a linear combination of $\ph_m$, $m\ge n+1$ with positive coefficients.
Now the equality
$$
I_0^\ga(\eta(\ph_n))=b_0^{\eta,\ga}S_{[1],n-1}(\ph_n)+c_0^{\eta,\ga}S_{[w_0],n-1}(\ph_n)
$$
implies $I_0^\ga(\eta(\ph_n))=(b_0^{\eta,\ga}-c_0^{\eta,\ga})S_{[1],n-1}(\ph_n)$ and thus
$$
|b_0^{\eta,\ga}-c_0^{\eta,\ga}|\ll\frac{|I_0^\ga(\eta(\ph_n))|}{n^{n-1}}.
$$
Similarly, we have $I_0^\ga(\eta(\ph_{n+1}))=(b_0^{\eta,\ga}+c_0^{\eta,\ga})S_{[1],n-1}(\ph_{n+1})$, and thus
$$
|b_0^{\eta,\ga}+c_0^{\eta,\ga}|\ll\frac{|I_0^\ga(\eta(\ph_n))|}{n^{n-1}}.
$$
So the same estimate holds for $b_0^{\eta,\ga}$ and $c_0^{\eta,\ga}$ separately.
We finish the proof of the lemma by estimating the period integral
$$
|I_0^\ga(\eta(\ph_n))|\ll n^4,
$$
independent of $\eta$.
To see this, let $\Delta$ denote the group Laplacian, i.e., $\Delta=-C+2C_K$, where $C$ is the Casimir-operator and $C_K$ is the Casimir operator of the group $K$.
Then for $\eta:\CD_{2n}\to L^2_\cusp(\Ga\bs G)$ a computation shows
$\Delta \eta(\ph_n)=n(n-1)\eta(\ph_n)$.
Since $\Delta$ is elliptic, positive definite and of order two on the 3-dimensional manifold $\Ga\bs G$, the operator $(1+\Delta)^{-2}$ has a continuous kernel $k(x,y)$.
It acts on the space of cusp forms through the continuous kernel
$$
k_0(x,y)=\sum_{j=1}^\infty (1+\la_j)^{-2}\phi_j(x)\ol{\phi_j(y)},
$$
where $\phi_j$ is an orthonormal basis of $L^2_\cusp(\Ga\bs G)$ consisting of $\Delta$-eigenfunctions and $\la_j$ is the eigenvalue of $\phi_j$.
Selberg has shown in the G\"ottingen lectures, that $k_0$ actually is an $L^2$-kernel.
We have
\begin{align*}
|\eta(\ph_n)(x)|&=(1+n(n-1))^2|(1+\Delta)^{-2}\eta(\ph_n)(x)|\\
&=(1+n(n-1))^2 \left|\int_{\Ga\bs G}k_0(x,y)\eta(\ph_n)(y)\,dy\right|\\
&\le (1+n(n-1))^2 \(\int_{\Ga\bs G}|k_0(x,y)|^2\,dy\)^{1/2}.
\end{align*}
As $k_0$ is an $L^2$-kernel, the latter integral is finite almost everywhere in $x$ and is locally bounded outside a set of measure zero.
Hence the continuous function $\eta(\ph_n)$ is locally bounded by a constant times $(1+n(n-1))^2$, so it is locally $O(n^4)$ and the same holds for the period integral $I_0^\ga(\eta(\ph_n))$.
\end{proof}

\section{Fourier expansion of Maa\ss\ forms}\label{sect3}
Let $f:\Ga\bs \H\to\C$ be a Maa\ss\ form, i.e., $f\in L^2(\Ga\bs \H)$ is an eigenform of the hyperbolic Laplacian, say $\Delta f=(\frac14-\la^2)f$ for $\la\in\C$.
Then there exists an automorphis representation $(\pi,V_\pi,\eta)$ with $\pi=\pi_\la$ such that $f=\eta(\ph_0)$, where $\ph_0\in V_\pi$ is a $K$-invariant function.
By scaling, one can achieve $\ph_0(K)=\{ 1\}$.
One then gets the \e{Fourier expansion}
$$
f(x)=\sum_{k\in\Z}a_k^{\eta,\ga} I_{\la,k\mu_\ga}^{\st}(\pi(\sigma^{-1} x)\ph_0).
$$
Note that the automorphic representation $\eta$, which determines the form $f$, enters on the right hand side only through the coefficients $(a_k^{\eta,\ga})_{k\in\Z}$.
In particular, if $x=\sigma a_t$ for $t\in\R$, 
one has $I_{\la,k\mu_\ga}^{\st}(\pi(\sigma^{-1} x)\ph_0)
=I_{\la,k\mu_\ga}^{\st}(\pi(a_t)\ph_0)=a_t^\mu I_{\la,k\mu_\ga}^{\st}(\ph_0)$, 
and therefore
$
f\(\sigma\smat{e^t}\ \ {e^{-t}}\)
$
equals
$$
\frac1{\Ga\(\frac{\la+1}2\)}
\sum_{k\in\Z}a_k^{\eta,\ga}e^{2\pi i\tilde\mu kt}
\Ga\(\frac{\la+2\pi ik\tilde\mu+1}4\)\Ga\(\frac{\la-2\pi ik\tilde\mu+1}4\).
$$

\section{Triple products}\label{sect4}
In this section we assume $\Ga$ to be cocompact. Then $L^2_\cusp(\Ga\bs G)=L^2(\Ga\bs G)$.
Let $(\pi,V_\pi,\eta)$ be an automorphic representation with $\pi=\pi_\la$ for some $\la\in\C$. 
Since $\pi$ is unitary, there is an anti-linear isomorphism to the dual $c: V_\pi\to V_\pi^*$.
Let $\breve\pi$ denote the dual representation on $V_\pi^*=V_{\breve\pi}$.
Let $\bar\cdot$ be the complex conjugation on $L^2(\Ga\bs G)$ and let $\breve \eta$ be the composition of the maps
$$
\xymatrix{
V_{\breve\pi}\ar[r]^{c^{-1}}&V_\pi\ar[r]^{\eta\ \ \ \ \ \ }&L^2(\Ga\bs G)\ar[r]^{\ol{\cdot}}&L^2(\Ga\bs G)
}
$$
Then $\breve \eta$ is a $G$-equivariant linear isometry of $V_{\breve\pi}$ into $L^2(\Ga\bs G)$, so $(\breve\pi,V_{\breve\pi},\breve \eta)$ is an automorphic representation as well.

Let $\Delta:\Ga\bs G\to\Ga\bs G\times\Ga\bs G$ be the diagonal map.
Let $\Delta^*:C^\infty(\Ga\bs G\times\Ga\bs G) \to C^\infty(\Ga\bs G)$ be the corresponding pullback map and
let $E=V_\pi^\infty\hat\otimes V_{\breve\pi}^\infty$, where $\hat\otimes$ denotes the projective completion of the algebraic tensor product.
Let 
$$
\eta_E: E\hookrightarrow C^\infty(\Ga\bs G)\hat\otimes C^\infty(\Ga\bs G)\cong C^\infty(\Ga\bs G\times\Ga\bs G)
$$ 
be given by $\eta\otimes\breve\eta$.
For $\ga$ as in the first section and $k\in\Z$
we get an induced functional on $E$,
$$
l_{\Delta(\ga)}^k= I_{k}^\ga\circ \Delta^*\circ \eta_E.
$$
In other words, for $w\in E$ we have
$$
l_{\Delta(\ga)}(w)=\frac 1{l(\ga)}\int_{A/\sp{a_\ga}}\eta_E(w)(\sigma_\ga a,\sigma_\ga a)\, da.
$$
This has the Fourier series expansion,
\begin{align*}
l_{\Delta(\ga)}(w) &= \sum_{k\in\Z} I_k^\ga\circ\eta\otimes I_{-k}^\ga\circ\breve\eta (w)\\
&= \sum_{k\in\Z}a_k^{\eta,\ga}a_{-k}^{\breve \eta,\ga}\,I_{\la,k\mu_\ga}^{\st}\otimes I_{-\la}^{-k\mu_\ga}(w)\\
&= \sum_{k\in\Z}|a_k^{\eta,\ga}|^2\,\hat w(k,-k),
\end{align*}
where the last line defines $\hat w$ and we have used the fact that $a_{-k}^{\breve \eta,\ga}=\ol{a_k^{\eta,\ga}}$.

Let $(\tau,V_\tau)$ be another element of $\hat G$. 
According to \cite{Mo} there is a canonical $G$-invariant continuous functional 
$$
T^\st_\tau : E\ \hat\otimes V_{\breve\tau}^\infty\ \to\ \C,
$$
and any other such functional is a scalar multiple of $T_\tau^\st$.
This induces a canonical $G$-equivariant continuous map
$$
T_\tau^\st : E\ \to\ V_\tau^\infty.
$$
On the other hand we have $\Delta^*\circ \eta_E\colon E\to L^2(\Ga\bs G)^\infty$.
For an automorphic representation $(\tau, V_\tau, \eta_\tau)$ we have an orthogonal projection $\Pr_{\eta_\tau}: L^2(\Ga\bs G)\to V_\tau$ and thus we get a map
$$
T_{\eta_\tau}^\aut= \Pr_{\eta_\tau}\circ \Delta^*\circ \eta_E
$$
from $E$ to $V_\tau^\infty$.
Hence there is a coefficient $c(\eta,\eta_\tau)\in\C$ such that
$$
T_{\eta_\tau}^\aut= c(\eta,\eta_\tau) T_\tau^\st.
$$
Fix a complete family $(\eta_j)$ of normalized, pairwise orthogonal automorphic representations $\eta_j:\pi_j\to L^2(\Ga\bs G)$.
Then the spectral expansion of $\Delta^*\circ \eta_E$ is
$$
\Delta^*\circ \eta_E= \sum_j c(\eta,\eta_j) T_{\pi_j}^\st.
$$
And hence, for $w\in E$,
\begin{align*}
l_{\Delta(\ga)}(w) &= \sum_{j:\pi_j\notin\hat G_{\rm ds}} c(\eta,\eta_j)a_{0}^{\eta_j,\ga} I_{\pi_j}^0(T_{\pi_j}^\st(w))\\
&+ \sum_{\stack{j:\pi_j\in\hat G_{\rm ds}}{\pi_j\cong\CD_{2n}}} c(\eta,\eta_j)b_{0}^{\eta_j,\ga}\( S_{[1],n-1}(T_{\pi_j}^\st(w))- S_{[w_0],n-1}(T_{\pi_j}^\st(w))\),
\end{align*}
where $\hat G_{\rm ds}\subset\hat G$ is the set of all discrete series representations of $G$.
So we conclude

\begin{lemma}\label{3.1}
\begin{multline*}
\sum_{k\in\Z}|a_k^{\eta,\ga}|^2\,\hat w(k,-k) = \sum_{j:\pi_j\notin\hat G_{\rm ds}} c(\eta,\eta_j)a_{0}^{\eta_j,\ga} I_{\pi_j}^0(T_{\pi_j}^\st(w))\\
+ \sum_{\stack{j:\pi_j\in\hat G_{\rm ds}}{\pi_j\cong\CD_{2n}}} c(\eta,\eta_j)b_{0}^{\eta_j,\ga} \(S_{[1],n-1}(T_{\pi_j}^\st(w))-S_{[w_0],n-1}(T_{\pi_j}^\st(w))\).
\end{multline*}
\end{lemma}

We make this a bit more explicit by plugging in special test functions.
The result is the following Theorem.

\begin{theorem}
[Summation formula]\label{thm4.2}
For $\al\in C_c^\infty(\R^\times)$ we have
$$
\sum_{k\in\Z}|a_k^{\eta,\ga}|^2
\hat\al(\tilde\la+k\tilde\mu_\ga)
=
\sum_{j:\pi_j\notin \hat G_{\rm ds}}c(\eta,\eta_j)a_0^{\eta_j,\ga}
\int_{\R^2}W_j(t,x)\al(\hat t_x)\,dt\,dx,
$$
where
\begin{align*}
W_j(t,x)&=\(\left|\frac{(e^{2t}+x)x}{(e^{2t}+x-1)(x-1)}\right|^{\frac12}+\left|\frac{(e^{2t}+x)x}{(e^{2t}+x-1)(x-1)}\right|^{-\frac12}\)^{\frac{\la_j-1}2}\\
&\times
\left|\frac{(e^{2t}+x)x}{(e^{2t}+x-1)(x-1)}\right|^{\frac12}
e^{(\la_j+1)t}
|x|^{\la-1}
|e^{2t}+x|^{-\la-1},\\
\hat t_x&=\frac 12\log\left|\frac{(e^{2t}+x)(x-1)}{(e^{2t}+x-1)x}\right|,\\
\end{align*}
\end{theorem}

\begin{proof}
Let $H_1=\smat 1\ \ {-1}$.
For $x\in\R$ write $n(x)=\smat 1x\ 1$.
For $\ph\in C_c^\infty(\R)$ we set
$$
f_\ph(w_0n(x))= \ph(-\sfrac12\log |x|).
$$
For given $\la\in\a^*$ the function $f_\ph$ extends uniquely to an element of $V_\la^\infty$.
For $a=\exp(tH_1)\in A$ we have $a^\la=e^{\la t}$. We define $\tilde\la=\frac{\la}{2\pi i}$, so that $a^\la=e^{2\pi i\tilde\la t}$.
We normalize the Haar measure $da$ on $A$ such that $\int_Ag(a)da=\int_\R g(\exp(tH_1)dt$.
Assuming that $n_0=\smat 1 1\ 1$, we compute for $k\in\Z$,
\begin{align*}
I_{\la,k\mu_\ga}^{\st}(f_\ph) &= \int_Af_\ph(w_0n_0^a) a^{-\la-1- k\mu_\ga}da\\
&= \int_\R\ph(t)\,e^{-t(2\pi i\tilde\la+1+2\pi i k\tilde\mu_\ga)}dt.
\end{align*}

Likewise, for $\phi\in C_c^\infty(\R^2)$ let $w_\phi\in E$ be defined by $w_\phi(w_0n(x),w_0n(y))=\phi(-\sfrac12\log |x|,-\sfrac12\log |y|)$.
Here we identify $E$ with $V_\la^\infty\hat\otimes V_{-\la}^\infty$.
We get
\begin{align*}
\hat w_\phi(k,-k)&= \int_{\R^2}\phi(t,s)\,e^{(s-t)(2\pi i\tilde\la+2\pi ik\tilde\mu_\ga)-(s+t)}dt\,ds
\end{align*}

Note that $w_\phi(w_0n(x),w_0n(y))$ vanishes in a neighborhood of $\{ xy=0\}$ as well as in a neighborhood of  $\{x=\infty\}\cup\{y=\infty\}$, which means that $w_\phi$ indeed lies  in $E$ and that
$$
S_{[1],n-1]}(T_{\pi_j}^\st(w_\phi))=S_{[w_0],n-1}(T_{\pi_j}^\st(w_\phi))=0
$$
for every $n\in\N$ and every $\pi_j\cong \CD_{2n}$.
So for these test functions ony the first sum on the right hand side of Lemma \ref{3.1} is present.

\begin{lemma}\label{lem4.3}
Let $g=\emat abcd$ and $g'=\emat{a'}{b'}{c'}{d'}$ be in $G$.
If $cc'dd'=0$, then $w_\phi(g,g')=0$.
Otherwise,
$$
w_\phi(g,g')= |c|^{-\la-1}|c'|^{\la-1} \phi\(\frac12(\log\left|\frac cd\right|,\frac12\log\left|\frac{c'}{d'}\right|\).
$$
\end{lemma}

\begin{proof}
Replacing $g$ by $\esmat\ {-1}1\ g$ or $g'$ by $\esmat\ {-1}1\ g'$ does not change $w_\phi(g,g')$.
We therefore can restrict our attention to the case $g,g'\in G^0$.
The Ansatz $g=anw_0n(x)$ with $a=\smat y\ \ {1/y}$ and $n=\smat 1z\ 1$ leads to $y=1/c$ and $x=d/c$ as well as $z=ac$.
This gives the claim.
\end{proof}

With $w=w_\phi$ as above, we want to compute $I_0^{\pi_j}(T_{\pi_j}^\st(w))$.
We first consider the case $\pi_j=\pi_{\la_j}$ for some $\la_j\in\C$.
Recall that the functional $T^\st:\pi\otimes\breve\pi\otimes\breve{\pi_j}\to\C$ maps a given $\ph=w\otimes f$ to
\begin{align*}
T^\st(\ph) &= \int_G\ph(w_0n_0y,w_0y,y)\, dy\\
&= \int_{ANK} w(w_0n_0ank,w_0ank)f(ank)\, da\,dn\,dk\\
&= \int_K\int_{AN} a^{-\la_j+1}w(w_0n_0ank,w_0 ank)da\,dn\,f(k)\,dk.
\end{align*}
The induced map $T_{\pi_j}^\st:E\to\pi_j$ is defined via the pairing
$$
\breve\pi_j\otimes\pi_j=\pi_{-\la_j}\otimes\pi_{\la_j}\ \to\ \C
$$
given by
$$
(f\otimes h)\ \mapsto\ \int_K f(k)h(k)\, dk.
$$
The resulting map $T_{\pi_j}^\st :\pi\otimes\breve\pi\to\pi_j$ therefore is given by
$$
T_{\pi_j}^\st(w)(k)=\int_{AN}a^{-\la_j+1}w(w_0n_0ank,w_0 ank)da\,dn.
$$
We have
$$
I_0^{\pi_j}(T_{\pi_j}^\st(w))= \int_A \ul a(w_0n_0a_1)^{\la_j+1} T_{\pi_j}^\st(w)(\ul k(w_0n_0a_1))\,da_1
$$
which equals
$$
\int_A\int_{AN}\ul a(w_0n_0a_1)^{\la_j+1} a^{-\la_j+1}
w(w_0n_0an\ul k(w_0n_0a_1),w_0an\ul k(w_0n_0a_1))\,da\,dn\,da_1.
$$
We write $\ul k(w_0n_0a_1)=\ul{an}(w_0n_0a_1)^{-1}w_0n_0a_1$ and use the change of variables $an\mapsto an\ul{an}(w_0n_0a_1)$ in the $AN$-integral.
For this we use the formula $\int_Hf(xy)\,dx=\Delta(y^{-1})\int_Hf(x)\,dx$ for Haar integration over the group $H=AN$ and the fact that the modular function $\Delta$ of $AN$ equals $\Delta(a_tn)=e^{-2t}$.
In this way we see that $I_0^{\pi_j}(T_{\pi_j}^\st(w))$ equals
$$
\int_A\int_{AN}\ul a(w_0n_0a_1)^{\la_j-1}a^{\la-\la_j}
w(w_0n_0anw_0n_0a_1,w_0nw_0n_0a_1)\,da\,dn\,da_1.
$$
Writing $\ol n_x=\smat 1\ {-x} 1=w_0n_xw_0$ this equals
$$
\int_A\int_{AN}\ul a(w_0n_0a_1)^{\la_j-1}a^{\la-\la_j}
w(\ol n_0a^{-1}\ol nn_0a_1,\ol nn_0a_1)\,da\,dn\,da_1.
$$
Writing $a=\esmat{e^t}\ \ {e^{-t}}$, $n=\esmat 1 x\ 1$, and $a_1=\esmat{e^s}\ \ {e^{-s}}$ we get
$$
\ol n_0a^{-1}\ol nn_0a_1=\emat {-e^{-t}}{-e^{-t}}{e^s(e^tx+e^{-t})}{e^{-s}(e^{t}(x-1)+e^{-t})}
$$
and
$$
\ol nn_0a_1=\emat {-e^{s}}{-e^{-s}}{e^sx}{e^{-s}(x-1)}.
$$
By Lemma \ref{lem4.3} we conclude
that $w_\phi(\ol n_0a^{-1}\ol nn_0a_1,\ol nn_0a_1)$ equals
$$
e^{-2s}
|e^tx+e^{-t}|^{-\la-1}|x|^{\la-1}
\phi\(s+\frac12\log\left|\frac{e^{-2t}+x}{e^{-2t}+x-1}\right|,s+\frac12\log\left|\frac x{x-1}\right|\).
$$
We now pick $\al,\beta\in C_c^\infty(\R)$ and set 
$$
\phi(x,y)=\al(x-y)\beta\(\frac{x+y}2\).
$$
With this choice, we get
$$
\hat w_\phi(k,-k)=\hat\al(\tilde\la+k\tilde\mu_\ga)\hat\beta\(\frac1{\pi i}\),
$$
and
$w_\phi(\ol n_0a^{-1}\ol nn_0a_1,\ol nn_0a_1)$ equals
\begin{align*}
&e^{-2s}
|e^tx+e^{-t}|^{-\la-1}|x|^{\la-1}
\al\(\frac 12\log\left|\frac{(e^{-2t}+x)(x-1)}{(e^{-2t}+x-1)x}\right|\)\\
&\times
\beta\(s+\frac14\log\left|\frac{(e^{-2t}+x)x}{(e^{-2t}+x-1)(x-1)}\right|\).
\end{align*}
We conclude that $I_0^{\pi_j}(T_{\pi_j}^\st(w))$ equals
\begin{align*}
&\int_{\R^3}(e^{2s}+e^{-2s})^{\frac{\la_j-1}2}e^{(\la-\la_j)t}e^{-2s}|e^tx+e^{-t}|^{-\la-1}|x|^{\la-1}\\
&\times
\al\(\frac 12\log\left|\frac{(e^{-2t}+x)(x-1)}{(e^{-2t}+x-1)x}\right|\)\\
&\times
\beta\(s+\frac14\log\left|\frac{(e^{-2t}+x)x}{(e^{-2t}+x-1)(x-1)}\right|\)\,dx\,ds\,dt.
\end{align*}
After the change of variables $s\mapsto s-\frac14\log\left|\frac{(e^{-2t}+x)x}{(e^{-2t}+x-1)(x-1)}\right|$
we end up with
\begin{align*}
&\int_{\R^3}\(e^{2s}\left|\frac{(e^{-2t}+x)x}{(e^{-2t}+x-1)(x-1)}\right|^{\frac12} +e^{-2s}\left|\frac{(e^{-2t}+x)x}{(e^{-2t}+x-1)(x-1)}\right|^{-\frac12}\)^{\frac{\la_j-1}2}\\
&\times
\left|\frac{(e^{-2t}+x)x}{(e^{-2t}+x-1)(x-1)}\right|^{\frac12} 
e^{(\la-\la_j)t} |x|^{\la-1}|e^tx+e^{-t}|^{-\la-1}\\
&\times
\al\(\frac 12\log\left|\frac{(e^{-2t}+x)(x-1)}{(e^{-2t}+x-1)x}\right|\)\beta(s)\,ds\,dt\,dx.
\end{align*}
Assuming $\beta\ge 0$ and $\int_\R\beta(s)\,ds=1$ we replace $\beta$ with $\beta_T(x)=T\beta(Tx)$ and let $T\to\infty$ to get the claim.
The interchange of limit and sum resp. integral is justified by the Theorem of dominated convergence as follows.
Using results from \cite{BR} one deduces that 
$$
c(\eta,\eta_j)=O(|\la_j|^{2+\eps}).
$$
Weyl's asymptotic law says that $|\la_j|=O(j)$ and with Theorem \ref{thm2.5} we get
$$
a_0^{\eta_j,\ga}c(\eta,\eta_j)=O|\la_j|^3)=O(j^3).
$$
So we need an estimate $O(|\la_j|^{-5})$ of the above integral, with $\beta=\beta_T$, which is uniform in $T$.
One achieves this by iterated use of the fact that $2\la_je^{2t\la_j}=\frac d{dt}e^{2t\la_j}$ and using integration by parts.
This, actually, is the place where it is needed that $\al$ be vanishing in a neighborhood of zero, as otherwise there would appear boundary terms of this integration by parts.
The Theorem follows.
\end{proof}

{\small Mathematisches Institut, 
Auf der Morgenstelle 10, 
72076 T\"ubingen, 
Germany, 
\tt deitmar@uni-tuebingen.de}

\today

\end{document}